\documentclass{amsart2000}
\usepackage{hyperref}

\theoremstyle{plain}
\newtheorem{thm}{Theorem}
\newtheorem{prop}[thm]{Proposition}

\newtheorem{cor}[thm]{Corollary}

\theoremstyle{remark}
\newtheorem{rem}[thm]{Remark}
\newtheorem*{ack}{Acknowledgments}

\renewcommand{\b}[1]{\mathbf{#1}}

\newcommand{\A}{{\mathcal A}}
\newcommand{\D}{{\mathcal D}}
\newcommand{\LL}{{\mathcal L}}
\newcommand{\R}{{\mathcal R}}
\newcommand{\C}{{\mathbb C}}
\newcommand{\ii}{{\sqrt{-1}}}
\newcommand{\T}{{({\mathbb C}^*)^n}}
\newcommand{\la}{{\lambda }}
\newcommand{\bl}{{\boldsymbol{\lambda}}}

\begin{document}

\title[Triples of arrangements and local systems]
{Triples of arrangements and local systems}

\author{Daniel C.~Cohen} 
\address{Department of Mathematics, Louisiana State University, Baton
Rouge, LA 70803}
\email{\href{mailto:cohen@math.lsu.edu}{cohen@math.lsu.edu}}
\urladdr{\href{http://www.math.lsu.edu/~cohen/}
{www.math.lsu.edu/\~{}cohen}}
\thanks{Partially supported by Louisiana Board of Regents grant
LEQSF(1999-2002)-RD-A-01 and by National Security Agency grant
MDA904-00-1-0038}



\begin{abstract} 
For a triple of complex hyperplane arrangements, there is a well-known
long exact sequence relating the cohomology of the complements.  We
observe that this result extends to certain local coefficient systems,
and use this extension to study the characteristic varieties of
arrangements.  We show that the first characteristic variety may
contain components that are translated by characters of any order,
thereby answering a question of A. Suciu.
\end{abstract}


\maketitle

\section{Introduction} \label{sec:intro}

Let $\A=\{H_{1},\dots,H_n\}$ be a hyperplane arrangement in $\C^\ell$,
with complement $M(\A)=\C^{\ell}\setminus\bigcup_{j=1}^{n}H_{j}$.  The
cohomology of $M(\A)$ with coefficients in a complex rank one local
system has been the subject of considerable recent interest, with
applications in the theory of multivariable hypergeometric functions
\cite{AK,Gel1,OT2} and mathematical physics \cite{Va}.  The study of
such local systems is also intimately related to a number of
algebraic, combinatorial, and topological invariants of arrangements.

Note that $M(\A)$ has the homotopy type of an $\ell$-dimensional cell
complex.  The set of rank one local systems on $M(\A)$ may be realized
as the complex torus $\T\cong H^1(M(\A);\C^*)$, with coordinates
$(z_1,\dots,z_n)$, $z_j$ corresponding to $H_j\in\A$.  A point
$\b{t}\in\T$ determines a representation $\rho:\pi_1(M(\A))\to\C^*$
given by $\gamma_j\mapsto t_j$ for any meridian loop $\gamma_j$ about
the hyperplane $H_j$ of $\A$, and a corresponding rank one local
system $\LL_{\b{t}}$ on $M(\A)$.  For sufficiently generic $\b{t}$,
the cohomology $H^{q}(M(\A);\LL_{\b{t}})$ vanishes, except possibly in
dimension $\ell$.  Those $\b{t}$ for which the cohomology does not
vanish comprise the characteristic varieties
\[ 
\Sigma^{q}_m(\A)=\{t\in\T \mid
\dim H^{q}(M(\A);\LL_{\b{t}})\ge m\}. 
\] 
These varieties are homotopy-type invariants of $M(\A)$, which inform
on the homology of certain abelian covers \cite{L1,Sa}, and also on
the structure of the second nilpotent quotient of the fundamental
group \cite{MS}.  It follows from work of Arapura~\cite{Ar} that the
characteristic varieties, or cohomology support loci, are unions of
torsion-translated subtori of $\T$ in more general circumstances.

There is a combinatorial cohomology theory, with corresponding support
loci, that is closely related to the study of (rank one) local systems
on complements of arrangements.  We use notation and results from
\cite{OT1}.  Let $A(\A)$ be the Orlik-Solomon algebra of $\A$
generated by the 1-dimensional classes $a_j$, $1\leq j\leq n$.  It is
the quotient of the exterior algebra generated by these classes by a
homogeneous ideal, so is a finite dimensional graded $\C$-algebra. 
There is an isomorphism of graded algebras $A(\A) \cong H^*(M;\C)$. 
In particular, $A^1(\A) \cong H^1(M;\C) \cong \C^n$, with
coordinates $(y_1,\dots,y_n)$, $y_j$ corresponding to $H_j\in\A$.  A
point $\bl \in \C^n$ determines an element $a_{\bl}=\sum_{j=1}^n \la_j
a_j$ in $A^{1}(\A)$, and multiplication by $a_{\bl}$ provides $A(\A)$
with the structure of a cochain complex.  For sufficiently generic
$\bl$, the cohomology $H^{q}(A(\A),a_{\bl}\wedge)$ vanishes, except
possibly in dimension $\ell$.  Those $\bl$ for which the cohomology
does not vanish comprise the resonance varieties
\[
\R^{q}_m(\A)=\{\bl\in\C^{n} \mid 
\dim H^{q}(A(\A),a_{\bl}\wedge)\ge m\}.  
\] 
These varieties are isomorphism-type invariants of $A(\A)$, which, if
$\A$ is central, inform on the structure (through rank 3) of the
matroid of $\A$ \cite{Fa}.  It follows from Theorem
\ref{thm:tangentcone} below that the resonance varieties are unions of
linear subspaces of $\C^{n}$.

Recent work of a number of authors \cite{CO,CScv,L2,L3,LY1}
establishes the following relationship between the characteristic and
resonance varieties of an arrangement.

\begin{thm} \label{thm:tangentcone}
Let $\A$ be an arrangement of $n$ hyperplanes in $\C^\ell$.  Then for
each $q$ and $m$, the resonance variety $\R_{m}^{q}(\A)$ coincides
with the tangent cone of the characteristic variety
$\Sigma_{m}^{q}(\A)$ at the point $\b{1}=(1,\dots,1)\in\T$.
\end{thm}
This result shows that the structure of the characteristic varieties
is largely combinatorial.  Components of $\Sigma_{m}^{q}(\A)$
containing $\b{1}$ are determined by the resonance varieties of the
Orlik-Solomon algebra $A(\A)$.  In turn, $A(\A)$ is determined by the
combinatorics of $\A$, namely the intersection poset $L(\A)$, see
\cite{OT1}.  The best understood characteristic and resonance
varieties are those for the first cohomology,
$\Sigma_m(\A):=\Sigma_m^1(\A)$ and $\R_m(\A):=\R_m^1(\A)$.  We refer
to \cite{CScv,Fa,L2,LY1} and the surveys \cite{F2,Yuz}.  In
particular, explicit combinatorial descriptions of the resonance
varieties $\R_m(\A)$, and hence of the components of $\Sigma_{m}(\A)$
containing $\b{1}$, are known.

However, as noted by Yuzvinsky \cite[\S8.4]{Yuz}, 
``relations between the resonance
and respective characteristic varieties are still mysterious,'' Theorem
\ref{thm:tangentcone} notwithstanding.  There are arrangements for which
the characteristic varieties contain components which do not pass through
$\b{1}$, so are not a priori detectable by combinatorial means.  The first
such example was found in \cite{CScv}.  The most striking examples of this
phenomenon were found recently by Suciu \cite{Su}, who exhibits
arrangements for which the first characteristic variety contains
positive-dimensional translated components.

This discovery prompted a number of questions in \cite{Su}.  One can
ask, for instance, are these translated components are combinatorially
determined?  what are the possible dimensions?  what are the possible
orders of translation?  In this note, we answer the last of these
questions as follows.

\begin{thm} \label{thm:tori}
For any positive integer $r\ge 2$, there is a hyperplane arrangement
$\A$ whose first characteristic variety $\Sigma_1(\A)$ contains
positive-dimensional components which are translated by characters of
order $r$.
\end{thm}

This result is established in Section \ref{sec:tori}.  Call a
component of $\Sigma_m(\A)$ essential if it is not contained in a
subtorus $z_j=1$, in the coordinates for $\T$ specified above. 
Analogously, call a component of $\R_m(\A)$ essential if it is not
contained in a hyperplane $y_j=0$, in the chosen coordinates for
$\C^n$.  We exhibit a family of arrangements $\{\D_r\}_{r \ge 2}$,
whose first characteristic varieties have essential components that
are translated as asserted.  The first member of this family, $\D_2$,
is the ``deleted $\rm{B}_3$ arrangement,'' the principal example found
by Suciu \cite{Su}.  Viewing the $\rm{B}_3$ Coxeter arrangement as the
reflection arrangement associated to the full monomial group
$G(2,1,3)$, one is led to study the characteristic varieties of
deletions $\D_r$ of the arrangements associated to the full monomial
groups $G(r,1,3)$ for any $r$.

For small $r$, we found translated components in the characteristic
varieties $\Sigma_1(\D_r)$ by direct calculation.  This note stems
from an attempt to find a conceptual, and perhaps combinatorial,
explanation for the existence of these translated components.  To this
end, we revisit a standard construction in arrangement theory, that of
a deletion-restriction triple, in the context of local systems.  In
Section~\ref{sec:triples}, we observe that a well-known exact sequence
relating the (constant coefficient) cohomology of the complements of
the arrangements in a triple extends to certain local systems.  This
observation facilitates the proof of Theorem \ref{thm:tori}.

\section{Triples} \label{sec:triples}

Let $\A$ be a hyperplane arrangement in $\C^{\ell}$.  The choice of a
distinguished hyperplane $H$ in $\A$ gives rise to a {\em triple} of
arrangements $(\A,\A',\A'')$, where $\A'=\A\setminus\{H\}$ and $\A''$
is the arrangement in $H\cong\C^{\ell-1}$ with hyperplanes $\{K\cap H
\mid K \in \A'\}$.

Let $M=M(\A)$, $M'=M(\A')$, and $M''=M(\A'')$ be the complements of
the arrangements in a triple.  These complements are related by
\[
M = M' \setminus (M'\cap H) \quad\text{and}\quad M'' = M' \cap H.
\]
The interrelations among the cohomology (rings) of $M$, $M'$, and $M''$
have been the subject of much study, see \cite{JR,OT1,Yuz}.  For
instance, the following is well-known.

\begin{prop}[{\cite[Cor.~5.81]{OT1}}] \label{prop:triple}
For a triple of complex hyperplane arrangements, there is a long exact
sequence in cohomology
\begin{equation*} 
\cdots \to H^{q}(M';\C) \to H^{q}(M;\C) \to H^{q-1}(M'';\C) 
\to H^{q+1}(M';\C) \to \cdots
\end{equation*}
\end{prop}

Less well-known is the fact that this result extends to certain
non-trivial local systems.  Let $i:M\to M'$ and $j:M''\to M'$ be the
natural inclusions and~let~$\LL'$ be a local system on $M'$.  Then there
is an induced local system $\LL=i^{*}\LL'$ on~$M$ with trivial monodromy
about the distinguished hyperplane $H$, and a local system
$\LL''=j^{*}\LL$ on $M''$ by restriction.  Call $(\LL,\LL',\LL'')$ a
triple of local systems on~$(M,M',M'')$.

\begin{thm} \label{thm:triple}
Let $(\A,\A',\A'')$ be a triple of arrangements, and let
$(\LL,\LL',\LL'')$ be a triple of complex rank one local systems on
$(M,M',M'')$.  Then there is a long exact sequence in local system
cohomology
\begin{equation} \label{eq:exseq}
\cdots \to H^{q}(M';\LL') \to H^{q}(M;\LL) \to H^{q-1}(M'';\LL'') 
\to H^{q+1}(M';\LL') \to \cdots
\end{equation}
\end{thm}

H. Terao informs us that this result was obtained 
independently by Y. Kawahara.

With the observation that $\LL$ has trivial monodromy about the
distinguished hyperplane $H$, the proof of Proposition~\ref{prop:triple}
in \cite{OT1} may be used to establish Theorem~\ref{thm:triple}. While
this result is apparently known, it does not appear that it has been
recorded in the literature, and least in this form.  So we provide a brief
proof.  More detailed arguments, in the case of trivial local
coefficients, may be found in \cite{JR,OT1,Yuz}.

\begin{proof}[Proof of Theorem \ref{thm:triple}]
Consider the long exact sequence of the pair $(M',M)$ with
coefficients in $\LL'$.  Since the local system $\LL'$ has trivial
monodromy about the distinguished hyperplane $H$, the restriction of
$\LL'$ to $M$ is itself a local system, $\LL'|_{M}=i^{*}\LL'=\LL$.  So
the long exact sequence of the pair is of the form
\[
\cdots \to H^{q}(M',M;\LL') \to H^{q}(M';\LL) \to H^{q}(M;\LL) 
\to H^{q+1}(M',M;\LL') \to \cdots.
\]
Thus it suffices to show that $H^{q+1}(M',M;\LL') \cong
H^{q-1}(M'';\LL'')$ for each $q$.

Let $E$ be a tubular neighborhood of $M''=M'\cap H$ in $M'$, where $H \in
\A$ is~the distinguished hyperplane. The neighborhood $E$ admits the
structure of a trivial bundle over $M''$ with fiber $\C$.  Identify the
zero section of this bundle with $M''$.  Then the complement of the zero
section, $E_{0}$, may be identified with $E \setminus M''$.

By excision, we have $H^{*}(M',M;\LL') \cong H^{*}(E,E_{0};\LL')$. 
Furthermore, the triviality of the bundle $E \to M''$ yields
$(E,E_{0}) \cong (M'' \times \C,M'' \times \C^{*}) \cong M'' \times
(\C,\C^{*})$.  The restriction, $\LL'|_{\C}$, of the local system
$\LL'$ to $(\C,\C^{*})$ is necessarily trivial.  Hence, we may compute
the cohomology $H^{*}(E,E_{0};\LL')$ using the K\"unneth formula:
\[
H^{q+1}(E,E_{0};\LL') \cong H^{q+1}(M'' \times (\C,\C^{*});\LL') 
\cong H^{q-1}(M'';\LL'|_{M''}) \otimes H^{2}(\C,\C^{*};\LL'|_{\C}).
\]
Since $H^{2}(\C,\C^{*};\LL'|_{\C})\cong \C$ and the restriction of
$\LL'$ to $M''$ is $\LL'|_{M''}=j^{*}\LL'=\LL''$, we have
$H^{q+1}(M',M;\LL')\cong H^{q-1}(M'';\LL'')$, as was required.
\end{proof}

\begin{cor} \label{cor:nonzero}
Let $(\LL,\LL',\LL'')$ be a triple of complex rank one local systems
on $(M,M',M'')$.  If $H^{q-1}(M'';\LL'')=0$ and $H^{q}(M;\LL)\neq 0$,
then $H^{q}(M';\LL')\neq 0$.  
If $\LL''$ is non-trivial, then $H^{1}(M';\LL')\cong H^{1}(M;\LL)$.
\end{cor}
\begin{proof} 
If $H^{q-1}(M'';\LL'')=0$, then the exact sequence
\eqref{eq:exseq} reduces to
\begin{equation} \label{eq:seq2}
\dots \to H^{q-2}(M'';\LL'')\to H^{q}(M';\LL')\to H^{q}(M;\LL)\to 0.
\end{equation}
If $H^{q}(M;\LL)\neq 0$, then
$H^{q}(M';\LL')\neq 0$ as well.  If $q=1$ and $\LL''$ is non-trivial, 
then $H^{i}(M'';\LL'')=0$ for $i < q$, so 
\eqref{eq:seq2} reduces further to
$H^{1}(M';\LL')\xrightarrow{\sim} H^{1}(M;\LL)$.
\end{proof}

\begin{rem} \label{rem:generalizations}
More general analogues of Theorem \ref{thm:triple} include the
following.  
\item[{(i)}] \quad Let $\LL'$ be a complex local system on $M'$ of rank
$m>1$, and let $\LL=i^{*}\LL'$ and $\LL''=j^{*}\LL'$ be the induced
local systems on $M$ and $M''$.  Then a straightforward modification
of the argument given above yields an exact sequence,
\[
\cdots \to H^{q}(M';\LL') \to H^{q}(M;\LL) \to H^{q-1}(M'';\LL'') 
\otimes \C^{m}
\to H^{q+1}(M';\LL') \to \cdots,
\]
in local system cohomology generalizing \eqref{eq:exseq}. 
\item[{(ii)}] [D. Massey]\quad More generally, let $X$ be a
(connected) complex analytic manifold.  Let $A$ be an closed analytic
subspace of $X$, and $H$ be a (connected) closed submanifold of $X$ of
(complex) codimension $c$ in $X$.  Denote by $i:X \setminus (A \cup H)
\to X \setminus A$ and $j:(X\setminus A) \cap H \to X \setminus A$ the
natural inclusions.

Let $F^\bullet$ be a complex of sheaves on $X\setminus A$ such that
the cohomology sheaves are locally constant.  Thus, the sheaves
$H^k(F^\bullet)$ may be non-zero in any (finite) number of degrees and
may have arbitrary ranks, but, for a fixed $k$, $H^k(F^\bullet)$ is a
local system on $X\setminus A$.

There is a fundamental distinguished triangle:
\[
j_!j^!F^\bullet\rightarrow F^\bullet\rightarrow 
Ri_*i^*F^\bullet\rightarrow j_!j^!F^\bullet[1].
\]
Since $j$ is the inclusion of one submanifold into another, and
$F^{\bullet}$ has locally constant cohomology, there is an isomorphism
$j^!F^\bullet\cong j^*F^\bullet[-2c]$.  Furthermore, for a closed
inclusion, $j_!\cong Rj_*$.  Thus, we obtain
\[
Rj_*j^*F^\bullet[-2c]\rightarrow F^\bullet\rightarrow 
Ri_*i^*F^\bullet\rightarrow (Rj_*j^*F^\bullet[-2c])[1].
\]
The associated long exact sequence in hypercohomology,
\[
\rightarrow H^{k-2c}(M'';F^\bullet_{|_{M''}})\rightarrow H^{k}(M';
F^\bullet)\rightarrow H^{k}(M; F^\bullet_{|_{M}}) \rightarrow
H^{k-2c+1}(M''; F^\bullet_{|_{M''}})\rightarrow,
\]
where $M=X\setminus (A\cup H)$, $M'=X \setminus A$, and
$M''=(X\setminus A) \cap H$, generalizes \eqref{eq:exseq}.
\end{rem}

\section{Translated Tori} \label{sec:tori}

We now use the results of the previous section to prove Theorem
\ref{thm:tori}.  For each positive integer $r\ge 2$, we exhibit an
arrangement whose first characteristic variety contains one or more
essential components translated by characters of order $r$.

Let $\A_r$ be the monomial arrangement of $n=3r+3$ hyperplanes in
$\C^3$ defined by
$Q(\A_r)=x_1^{}x_2^{}x_3^{}(x_1^r-x_2^r)(x_1^r-x_3^r)(x_2^r-x_3^r)$. 
These are the reflection arrangements associated to the full monomial
groups $G(r,1,3)$, see \cite{OT1}.  The characteristic varieties of
these arrangements were studied in \cite[\S6]{CScv}.  There, it is
shown that $\Sigma_{1}(\A_r)$ has an essential two-dimensional
component $C$ which contains the identity $\b{1}$.  Let
$\zeta=\exp(2\pi\ii/r)$.  Denote the coordinates of the complex torus
$\T$ by $(z_1,z_2,z_3,z_{12:1},\dots,z_{12:r},z_{13:1},\dots,
z_{13:r},z_{23:1},\dots, z_{23:r})$, where $z_i$ corresponds to the
hyperplane $H_i=\ker(x_i)$ and $z_{ij:k}$ to $H_{ij:k}=\ker(x_i -
\zeta^k x_j)$.  In these coordinates, the component $C\subset
\Sigma_{1}(\A_r)$ is given by
\[
C=\bigl\{\left(u^r,v^r,w^r,w,\dots,w,v,\dots,v,u,\dots,u\right) \in 
\T \mid uvw=1\bigr\}.
\]
Notice 
that setting, say, $w$ equal to $\zeta^q$ defines a
one-dimensional subvariety of $C$ which does not contain the 
identity, and that $z_3=1$ at every point in this subvariety.

In light of this observation, let $\D_r$ be the ``monomial deletion''
$\A_r \setminus \{H_3\}$, with defining polynomial $Q(\D_r)= x_1 x_2
(x_1^r - x_2^r)(x_1^r - x_3^r)(x_2^r - x_3^r)$.  Note that
$|\D_r|=n-1=3r+2$, and denote the coordinates of $(\C^*)^{n-1}$ by
$z_1$, $z_2$, and $z_{ij:k}$, ordered as above.  We assert that the
subvariety $C_q$ of $(\C^*)^{n-1}$ defined by
\begin{equation*} \label{eq:translatedtorus}
C_q=\bigl\{\left(u^r,v^r,\zeta^q,\dots,\zeta^q,v,
\dots,v,u,\dots,u\right)
\in (\C^*)^{n-1} \mid uv\zeta^q=1\bigr\}
\end{equation*}
is a component of the first characteristic variety $\Sigma_{1}(\D_r)$
for each $q$, $1\le q\le r-1$, that is clearly essential.  Let $T_q$
be the one-dimensional subtorus of $(\C^*)^{n-1}$ given~by
\[
T_q=\bigl\{\left(u^r,u^{-r},1,\dots,1,u^{-1},
\dots,u^{-1},u,\dots,u\right)
\in (\C^*)^{n-1} \mid u \in\C^*\bigr\},
\]
and let $\tau_{q} =
(1,1,\zeta^q,\dots,\zeta^q,\zeta^{-q},\dots,\zeta^{-q},1,\dots,1)$. 
Checking that $C_q = \tau_{q} \cdot T_q$, we see that $C_q$ is a
translated subtorus of $(\C^*)^{n-1}$, translated by a character of
order~$r$.

We first establish the containment $C_q \subset \Sigma_{1}(\D_r)$. 
For this, let $\b{t}'$ be a point in $C_q$, given by
$\b{t}'=(t_1',t_2',t_{12:1}',\dots,t_{12:r}',t_{13:1}',\dots,
t_{13:r}',t_{23:1}',\dots,t_{23:r}'\bigr)$, and let $\LL'$ denote the
corresponding rank one local system on $M(\D_r)$.  The induced local
system $\LL=i^*\LL'$ on $M(\A_r)$, where $i:M(\A_r) \to M(\D_r)$ is
the inclusion, corresponds to the point $\b{t} \in \T$ given by
$t_1=t_1'$, $t_2=t_2'$, $t_3=1$, and $t_{ij:k}=t_{ij:k}'$ for all $i$,
$j$, and $k$.  Since $\b{t}' \in C_q$, we have $\b{t} \in C \subset
\Sigma_{1}(\A_r)$.  Thus, $H^1(M(\A_r);\LL) \neq 0$.

Let $\A_r''=\{H_3 \cap H \mid H \in \D_r\}$ be the restriction
corresponding to the deletion $\D_r = \A_r \setminus \{H_3\}$.  By
Corollary \ref{cor:nonzero}, to show that $H^1(M(\D_r);\LL') \neq 0$,
i.e., that $\b{t}' \in \Sigma_{1}(\D_r)$, it suffices to show that the
local system $\LL'' = j^*\LL'$ on $M(\A_r'')$ is non-trivial, where
$j:M(\A_r'') \to M(\D_r)$ is the inclusion.  The hyperplanes of
$\A_r''\subset H_3$ are $H_1'' = H_3 \cap H_1 \cap H_{13:1} \cap \dots
\cap H_{13:r}$, $H_2'' = H_3 \cap H_2 \cap H_{23:1} \cap \dots \cap
H_{23:r}$, and $H_{12:k}'' = H_3 \cap H_{12:k}$ for $1\le k \le r$. 
If $\b{t}''=(t_1'',t_2'',t_{12:1}'',\dots,t_{12:r}'') \in
(\C^*)^{r+2}$ records the monodromy of $\LL''$ about these
hyperplanes, then $t_1'' = t_3 t_1 t_{13:1} \cdots t_{13:r} = 1$,
$t_2'' = t_3 t_2 t_{23:1} \cdots t_{23:r} = 1$, and $t_{12:k}'' = t_3
t_{12:k} = \zeta^q$ for $1 \le k\le r$.  Thus, $\LL''$ is non-trivial,
and $H^1(M(\D_r);\LL') \neq 0$.  Hence, $\b{t}' \in \Sigma_{1}(\D_r)$,
and $C_q \subset \Sigma_{1}(\D_r)$.

It remains to show that $C_q$ is a translated component of
$\Sigma_{1}(\D_r)$, i.e., that $C_q$ is not contained in a component
with passes through the identity $\b{1} \in (\C^*)^{n-1}$.  Suppose
otherwise.  Then $C_q \subset S \subset \Sigma_{1}(\D_r)$, and $\b{1}
\in S$.  Since $C_{q}$ is essential, so is $S$.  There is then an
essential component $L \subset \R_{1}(\D_{r})$ corresponding to $S$ by
Theorem \ref{thm:tangentcone}.  However, the characterizations of the
resonance varieties in \cite{Fa,LY1} rule out the existence of such a
component $L$.  For instance, it is readily checked that (the matroid
of) the arrangement $\D_r$ admits no non-trivial neighborly partition. 
It follows that all components of $\R_1(\D_r)$ are non-essential,
cf.~\cite[\S3]{Fa}.  Thus, $C_q$ is indeed a translated component of
$\Sigma_{1}(\D_r)$, which completes the proof of
Theorem~\ref{thm:tori}.

\smallskip

We conclude by noting a consequence of the above discussion.

\begin{cor} 
For any positive integer $n$, there is a hyperplane arrangement $\A$
for which $\Sigma_{1}(\A)$ contains $n$ essential positive-dimensional
translated tori.
\end{cor}

\begin{ack}
We thank D.~Massey for communicating to us the generalization of
Theorem~\ref{thm:triple} found in Remark \ref{rem:generalizations}
(ii).  We also thank M. Falk, P. Orlik, A. Suciu, and H. Terao for
useful conversations.
\end{ack}

\end{document}